\documentclass[12pt]{article}
\usepackage{url}
\usepackage{amsmath, amssymb}

\usepackage{color}
\usepackage{multirow}

\usepackage[margin=25mm]{geometry}

\newcommand{\Gd}{{\mathbb Z}^2_{\diamond}}

\renewcommand{\c}{c}
\newcommand{\Z}{\mathbb{Z}}

\newcommand{\Ztaxi}{\vec{\mathbb Z}^2}

\renewcommand{\Gd}{{\mathbb Z}^2_{\Diamond}}

\newcommand{\qed}[0]{{\hspace*{\fill}\mbox{$\Box$}}}

\newcommand{\besttorus}{5.3506}
\newcommand{\bestbox}{7.0852} 
\newcommand{\bestmuupper}{1.58746}
\newcommand{\bestmulower}{1.55701}   
\newcommand{\bestlambdalower}{4.8771}        

\newtheorem{theorem}{Theorem}[section]
\newtheorem{definition}{Definition}[section]
\newtheorem{lemma}[theorem]{Lemma}

\makeatletter
\def\bigpar{\bigbreak\@afterindentfalse\@afterheading\ignorespaces}
\def\medpar{\medbreak\@afterindentfalse\@afterheading\ignorespaces}
\def\smallpar{\smallbreak\@afterindentfalse\@afterheading\ignorespaces}

\makeatletter
\def\eqalign#1{\,\vcenter{\openup\jot\m@th
  \ialign{\strut\hfil$\displaystyle{##}$&$\displaystyle{{}##}$\hfil
     \crcr#1\crcr}}\,}
\makeatother

\begin{document}

\title{Phase Coexistence for the Hard-Core Model on $\Z^2$}
\author{Antonio Blanca\thanks{School of Computer Science, Georgia Institute of Technology, Atlanta GA 30332; ablanca@cc.gatech.edu. Research supported in part by the NSF grants CCF-1420934, CCF-1555579 and CCF-1617306.}, ~Yuxuan Chen\thanks{Computer Science Department, Columbia University, New York NY 10027; yuxuan.chen@columbia.edu. Research supported in part by the
Simons Foundation grant 360240.}, \\David Galvin\thanks{Department of Mathematics, University of Notre Dame, Notre Dame IN 46656; dgalvin1@nd.edu. Research supported in part by the
Simons Foundation grant 360240 and by the National Security Agency grant NSA H98230-13-1-0248.}, ~Dana Randall\thanks{School of Computer Science, Georgia Institute of Technology, Atlanta GA 30332; randall@cc.gatech.edu. Research supported in part by the NSF grant CCF-1526900.} ~and Prasad Tetali\thanks{School of Mathematics, Georgia Institute of Technology, Atlanta, GA 30332; tetali@math.gatech.edu. Research supported in part by the NSF grant DMS-1407657.}}
\date{\today}

\maketitle

\begin{abstract}
The hard-core model has attracted much attention across several disciplines, representing lattice gases in statistical physics and independent sets in discrete mathematics and computer science. On finite graphs, we are given a parameter $\lambda$, and an independent set $I$ arises with probability proportional to $\lambda^{|I|}$.  On infinite graphs a Gibbs measure is defined as a suitable limit with the correct conditional probabilities, and we are interested in determining when this limit is unique and when there is phase coexistence, i.e., existence of multiple Gibbs measures.

It has long been conjectured that on $\Z^2$ this model has a critical value $\lambda_c \approx 3.796$ with the property that if $\lambda < \lambda_c$ then it exhibits uniqueness of phase, while if $\lambda > \lambda_c$ then there is phase coexistence. 
Much of the work to date on this problem has focused on the regime of uniqueness, with the state of the art being recent work of Sinclair,  Srivastava, \v{S}tefankovi\v{c} and Yin showing that there is a unique Gibbs measure for all $\lambda < 2.538$. Here 
we explore the other direction and prove that there are multiple Gibbs measures for all $\lambda > \besttorus$. We also show that with the methods we are using we cannot hope to replace $\besttorus$ with anything below $\bestlambdalower$.

Our proof begins along the lines of the standard Peierls argument, but we add two innovations. First, following ideas of Koteck\'y and Randall, we construct an event that distinguishes two boundary conditions and always has long contours associated with it, obviating the need to accurately enumerate short contours. Second, we obtain improved bounds on the number of contours by relating them to a new class of self-avoiding walks on an oriented version of $\Z^2$.

\bigskip\noindent
2010 \textit{Mathematics subject classification:} 60C05, 68R05

%
\end{abstract}

\thispagestyle{empty}

\newpage

\setcounter{page}{1}

\section{Introduction}

For a graph $G$ let ${\mathcal I}(G)$ denote the set of independent sets of $G$. For finite $G$ the {\em hard-core measure on $G$ with parameter $\lambda$} is the measure $\mu_{G,\lambda}$ supported on  ${\mathcal I}(G)$ given by $\mu_{G,\lambda}(I) \propto \lambda^{|I|}$ for each $I \in {\mathcal I}(G)$, or equivalently
$$
\mu_{G,\lambda}(I) = \frac{\lambda^{|I|}}{\sum_{J \in {\mathcal I}(G)} \lambda^{|J|}}.
$$
The hard-core measure is a simple mathematical model of a gas with particles
of non-negligible size. The vertices of $G$ are regarded
as positions, each of which can be occupied by a particle, subject to the
rule that two neighboring sites cannot both be occupied (particles cannot
overlap).

On infinite graphs, 
which may admit infinitely many independent sets, we make sense of the notion of choosing an independent set $I$ with probability proportional to $\lambda^{|I|}$ using the machinery of Gibbs measures. Roughly speaking, these are measures supported on the set of independent sets of $G$ whose conditional restrictions to finite subgraphs agree with the (suitably conditioned) finite hard-core measure.


Formally, let $G=(V,E)$ be infinite and locally finite  (i.e.\ no vertices of infinite degree). 
We say that a property holds for \textit{$\mu$-almost every} independent set in $\mathcal{I}(G)$, if the set $\mathcal{A}$ of independent sets for which the property does not hold has measure $0$; i.e., $\mu(\mathcal{A}) = 0$. 
\begin{definition} \label{def-Gibbs}
	A probability measure $\mu$ is a {\em Gibbs measure} for the hard-core model with parameter $\lambda$ on $G$ if,
	for every finite $\Lambda \subset V$, every $J \in  {\mathcal I}(G)$ and $\mu$-almost every $I \in  {\mathcal I}(G)$,
	 \[
	 \mu( J \mid I \cap \bar{\Lambda}) = \left.
	 \begin{cases}
	 \frac{1}{Z_{\Lambda,\lambda}^I} \lambda^{|J \cap \Lambda|}& \textrm{if}~\,I \cap \bar{\Lambda} = J \cap \bar{\Lambda};  \\ 
	 \\
	  0         & \textrm{otherwise},
	 \end{cases}
	 \right.
	 \]
 	 where $\bar{\Lambda} = V \setminus \Lambda$ and
 	 \[
 	 Z_{\Lambda,\lambda}^I = \sum_{K \in {\mathcal I}(G):\,K \cap \bar{\Lambda} = I \cap \bar{\Lambda}}   \lambda^{|K \cap \Lambda|}.
 	 \]
\end{definition}

\noindent
See, for example, \cite{Georgii,Weitz} for a very thorough treatment of this topic.

General compactness arguments show that an infinite, locally finite 
graph $G$ admits at least one Gibbs measure. A central concern of statistical physics
(again see \cite{Georgii} for a thorough discussion) is understanding when a particular system --- the independent set model in the present setting --- exhibits {\em phase coexistence}
(also known as {\em phase transition}) on a given
infinite $G$, meaning that it admits more than one
Gibbs measure.

The canonical (and by far most studied, and physically relevant) case of the hard core measure is
that of the usual nearest neighbor graph on the integer grid $\Z^d$. In this note we specifically consider the two-dimensional grid $\Z^2$. Formally this is the graph whose vertex set is the set of pairs $(x,y) \in \Z \times \Z$, with two pairs adjacent if they differ on exactly one coordinate, and differ by $\pm 1$ on that coordinate. 

For the classical Ising model, seminal work of Onsager \cite{onsager} established the precise value ($\beta_c(\Z^2)=\log(1+\sqrt{2})$) of the critical inverse temperature below which that model exhibits uniqueness of phase and above which it exhibits phase coexistence.  Only recently have the analogous values for the (more general) $q$-state Potts model been established, in work of Beffara and Duminil-Copin \cite{bd-c}, settling a more than half-a-century old open problem. 

Such precise results for the hard-core model seem far out of reach with currently available methods. It has long been conjectured, though --- with computational support, see e.g. \cite{BaxterEntingTsang} --- that there is a $\lambda_c \approx 3.796$ such that the hard-core model on $\Z^2$ exhibits phase coexistence for all $\lambda > \lambda_c$, but not for any $\lambda < \lambda_c$.

Starting with Dobrushin \cite{dob} in 1968, both physicists and mathematicians have been developing techniques to
approach this conjecture. Most of the attention has
focused on establishing ever larger values of $\lambda$ below which there is uniqueness
of phase. The problem has proved to be a fruitful one for the blending of ideas from
physics, discrete probability and theoretical computer science, with improvements to
our understanding having been made successively by Radulescu and
Styer \cite{radulescustyer}, van den Berg and Steif \cite{bs}, Weitz \cite{weitz}, Restrepo, Shin, Tetali, Vigoda and Yang \cite{rstvy}, and Vera, Vigoda and Yang \cite{VeraVigodaYang}, among others. The state of the
art is recent work of Sinclair, Srivastava, \v{S}tefankovi\v{c} and Yin \cite{sssy}, building on the novel ideas of Weitz, which establishes that there is a unique Gibbs measure for all $\lambda < 2.538$.

Much less is known about the regime of phase coexistence. Dobrushin \cite{dob} established
that there is a $C>0$ such that for all $\lambda > C$ there are multiple Gibbs measures. He did not explicitly calculate $C$, but around the time
of the writing of \cite{bc7} Borgs \cite{borgs-pc} reported that 
a direct implementation of Dobrusin's argument would yield a value of $C$ strictly greater than $80$.
A later computation by the third author showed that
$C \approx 300$. 

Our main aim in this paper is to give a reasonable upper bound on $\lambda_c$ (if it exists).
\begin{theorem} \label{thm-phasecoexistence}
The hard-core model on $\Z^2$ with activity $\lambda$ admits multiple Gibbs measures for all $\lambda > \besttorus$.
\end{theorem}

An intuition for the meaning of multiple Gibbs measures can be gleaned from the following recipe for producing them. For an independent set $I$ and a finite set $W \subseteq V$, let ${\mathcal I}^I(W)$ be the (also finite) set of independent sets that agree with $I$ off $W$.
Fix $I$ and a nested sequence $(W_i)_{i=1}^\infty$ of
finite subsets of $V$ satisfying $\cup_i W_i=V$. For each $i$ let $\mu^I_i$
be the measure on ${\mathcal I}^I(W_i)$ in which each $J$ is selected with probability proportional to $\lambda^{|J\cap W_i|}$. Any
(weak) subsequential limit of the $\mu^I_i$'s (and by compactness there
must be at least one such) is a Gibbs measure. This fact was originally
proved, in a much more general context, by Dobrushin \cite{Dobrushin4}; see
for example \cite{bs} for a treatment specific to the hard-core model on the lattice, or \cite[Theorem 3.5]{BrightwellWinkler} for a simple proof in the slightly more general context of graph homomorphism models. 

From this recipe we see that an interpretation of the existence of multiple Gibbs measures on $\Z^2$ is that the local behavior of a randomly chosen independent set in a box can be made to depend on a boundary condition imposed on the box, even in the limit as the size of the box grows to infinity.  This leads to what turns out to be the standard approach to showing multiple Gibbs measures, which is to consider the limiting measures corresponding to two different boundary conditions on boxes in the lattice centered at the origin, and to find a statistic that separates these two limits. For the hard-core model, it suffices (see \cite{bs}) to compare the even boundary condition --- all vertices on the boundary of a box at an even distance from the origin are occupied --- and its counterpart odd boundary condition, and the distinguishing statistic is typically the occupation of the origin. Under odd boundary condition the origin should be unlikely to be occupied, since independent sets with odd boundary and (even) origin occupied must have a contour --- a two-layer thick unoccupied ring of vertices separating an inner region around the origin that is in ``even phase'' from an outer region near the boundary that is in ``odd phase''. For large enough $\lambda$, such an unoccupied layer is costly, and so such configurations are unlikely. This is essentially the Peierls argument for phase coexistence, and was the approach taken by Dobrushin \cite{dob}.

As we will see presently, the effectiveness of the Peierls argument is driven by the number of contours of each possible length -- better upper bounds on the number of contours translate directly to better upper bounds on $\lambda_c$. Previous work of other researchers on phase coexistence in the hard-core model on $\Z^2$ had viewed contours as simple polygons in $\Z^2$, which are closely related to the very well studied family of self-avoiding walks. While this is essentially the best possible point of view when applying the Peierls argument on the Ising model, it is far from optimal for the hard-core model. The first contribution of the present paper is the realization that hard-core contours, if appropriately defined, can be viewed as simple polygons in the oriented Manhattan lattice, in which edges of $\Z^2$ that are parallel to the $x$-axis (respectively, the $y$-axis) are oriented positively if their $y$-coordinate (respectively, $x$-coordinate) is even, and negatively otherwise, with the additional constraint that contours cannot make two consecutive turns. The number of such polygons can be understood by analyzing a new class of self-avoiding walks, that we refer to as {\em taxi walks}. 
The number of taxi walks turns out to be significantly smaller than the number of ordinary self-avoiding walks, leading to significantly better bounds on $\lambda_c$.

There is a single number $\mu_{\rm taxi}>0$, the taxi walk connective constant, that asymptotically controls the number $\c_n$ of taxi walks of length $n$, in the sense that $\c_n = \mu_{\rm taxi}^n f_{\rm taxi}(n)$ with $f_{\rm taxi}(n)$ sub-exponential. Adapting methods of Goulden and Jackson \cite{GJ} we obtain 
estimates on $\mu_{\rm taxi}$, giving us good understanding of $\c_n$ for large $n$. The sub-exponential correction makes it difficult to control $\c_n$ for small $n$, however, presenting a major stumbling block to the effectiveness of the Peierls argument as we have just described it. Using the statistic ``occupation of origin'' to distinguish the two boundary conditions, one inevitably has to control $\c_n$ for both small and large $n$. The lack of precise information about the number of short contours leads to discrepancies between asymptotic
and actual bounds, such as that between the lower bound $C>80$ and $C\approx 300$ from a direct implementation of Dobrushin's argument for phase coexistence on $\Z^2$ discussed earlier. 

The second contribution of the present paper is the idea of using an event to distinguish the two boundary conditions that has the property that every independent set in the event has associated with it a long contour. This allows us to focus exclusively on the asymptotic growth rate of contours/taxi walks, and obviates the need for an analysis of short contours. 

In an earlier version of this work \cite{BGRT} we used the technology of fault lines, introduced in \cite{randall} to define a distinguishing event. After a talk by one of the authors, Koteck\'y pointed out an alternate approach \cite{kot-pc}. Consider a box $B$ of fixed size centered at the origin, and say that an independent set is even (respectively, odd) on the box if every vertex in $B$ at an even (respectively, odd) distance from the origin is either in the independent set, or potentially could be in the sense that none of its neighbors are. Then consider the event that a randomly chosen independent set drawn from a much larger box (with boundary condition) is even on $B$, conditioned on the event that it is either even or odd on $B$. Running the Peierls argument on this event leads to contours that completely encircle $B$, and so can be made arbitrarily long by choosing $B$ to be sufficiently large.        

In Section \ref{sec-taxi-walks} we introduce the notion of taxi walks and the taxi walk connective constant, that will be key to the precise theorem we prove (Theorem \ref{thm-phasecoexistence-taxi}), from which Theorem \ref{thm-phasecoexistence} follows via some numerical computation.  The proof of  Theorem \ref{thm-phasecoexistence-taxi} is then given in Section \ref{sec-proof_phase_coexistence}. In Section \ref{sec-taxibounds} we give the details of our upper and lower bounds on the taxi walk connective constant, and we conclude in Section \ref{sec-remarks} with some remarks.

\section{Taxi walks} \label{sec-taxi-walks}

Let $\Ztaxi$ be an orientation of $\Z^2$ in which an edge parallel to the $x$-axis (respectively, $y$-axis) is oriented in the positive $x$-direction if its $y$-coordinate is even (respectively, oriented in the positive $y$-direction if its $x$-coordinate is even), and is oriented in the negative direction otherwise. 
It is common to refer to $\Ztaxi$ as the {\em Manhattan lattice}: streets are horizontal, with even numbered streets oriented west to east and odd numbered streets oriented east to west, and avenues are vertical, with even numbered avenues oriented south to north and odd numbered avenues oriented north to south. 
\begin{definition}
A {\em self-avoiding walk} in $\Ztaxi$ of {\it length} $n$ {\it starting} at vertex $v_0$ is a sequence $v_0, v_1, \ldots, v_n$ of distinct vertices with, for each $i=1, \ldots, n$, $v_{i-1}v_i$ an edge of $\Ztaxi$ oriented from $v_{i-1}$ to $v_i$. The walk {\em turns} at $v_i$ ($1 \leq i \leq n-1$) if edges $v_{i-1}v_i$ and $v_iv_{i+1}$ are perpendicular, and {\em goes straight} if these edges are parallel.  A {\em taxi walk} is a self-avoiding walk in $\Ztaxi$ that does not turn at two consecutive vertices.
\end{definition}
We call these taxi walks because a savvy passenger in a Manhattan cab would be suspicious if the cab took two consecutive turns.

Let $\c_n$ be the number of taxi walks of length $n$ starting at the origin. A critical step in our arguments will be bounding $\c_n$. An easy upper bound is $\c_n \leq 2^n$, since there are always at most two ways to extend a taxi walk of length $n-1$, and an easy lower bound is $2^{n/2} \leq \c_n$ (for even $n$) and $2^{(n+1)/2} \leq \c_n$ (for odd $n$), since walks that always take two steps north at a time or two steps east at a time can always be extended in exactly two ways.
With very little extra work we can get a significantly better upper bound.
\begin{lemma}  \label{lem-fib}
$\c_n = O\left(((1 + \sqrt{5})/2)^n\right)$. 
\end{lemma}

\noindent {\bf Proof}:
A taxi walk $v_0, v_1, \ldots, v_n$ (with $v_0$ the origin) can be encoded by a pair $(a,\sigma)$, where $a \in \{N,E\}$ and $\sigma$ is a sequence of length $n-1$ over alphabet $\{s,t\}$, as follows: if $v_1=(0,1)$ then $a=N$ and if $v_1=(1,0)$ then $a=E$, and if the walk goes straight at $v_i$ then the $i$th entry of $\sigma$ is $s$, whereas it is $t$ if the walk turns at $v_i$. Distinct taxi walks evidently get distinct codes. It is well know that the number of sequences of length $n-1$ over alphabet $\{s,t\}$ without consecutive occurrences of the character $t$ is the $(n+1)$st Fibonacci number $f_{n+1}$ (defined by $f_0=0$, $f_1=1$, $f_n=f_{n-1}+f_{n-2}$ for $n \geq 2$). It follows that $\c_n \leq 2f_{n+1}=O\left(((1 + \sqrt{5})/2)^n\right)$.
\qed

\medskip

Using more sophisticated tools we can improve our bounds. In what follows we say that a function $f(n)$ defined on positive integers grows {\em sub-exponentially} if for all $\varepsilon > 0$ there is $n(\varepsilon)$ such that for all $n > n(\varepsilon)$ we have $f(n) < (1+\varepsilon)^n$. 
\begin{theorem} \label{thm-best_mut}
There is a constant $\mu_{\rm taxi}$ (the {\em taxi walk connective constant}) with $\bestmulower < \mu_{\rm taxi} < \bestmuupper$ and $\bestlambdalower < \mu_{\rm taxi}^4-1 < \besttorus$, and a function $f_{\rm taxi}(n)$ that grows sub-exponentially, such that $\c_n = f_{\rm taxi}(n) \mu_{\rm taxi}^n$ for all $n$. 
\end{theorem}
We defer the proof of Theorem \ref{thm-best_mut} to Section \ref{sec-taxibounds}. We end this section with a precise statement of our main theorem, from which Theorem \ref{thm-phasecoexistence} follows via Theorem \ref{thm-best_mut}.
\begin{theorem} \label{thm-phasecoexistence-taxi}
The hard-core model on $\Z^2$ with activity $\lambda$ admits multiple Gibbs measures for all $\lambda > \mu_{\rm taxi}^4 -1$.
\end{theorem}
Our bounds on $\mu_{\rm taxi}^4 -1$ involve the theory of irreducible bridges and the Goulden--Jackson cluster method, as well as extensive computation, but the proof of the existence of $\mu_{\rm taxi}$ and $f_{\rm taxi}(n)$ at the beginning of Section \ref{sec-taxibounds} and the bounds $\sqrt{2} \leq \mu_{\rm taxi} \leq (1+\sqrt{5})/2$ in this section are straightforward. A consequence of this is that if we wish to avoid using irreducible bridges, the Goulden--Jackson cluster method and computer-aided computations, we have via Theorem \ref{thm-phasecoexistence-taxi} and Lemma \ref{lem-fib} a weaker version of Theorem \ref{thm-phasecoexistence} that is still significantly better than any previous result; namely, that the hard-core model on $\Z^2$ with activity $\lambda$ admits multiple Gibbs measures for all $\lambda > (5+3\sqrt{5})/2$, and so for all $\lambda > 5.8542$.

\medskip

A consequence of our lower bound on $\mu_{\rm taxi}$ is that our present approach to phase coexistence cannot prove anything better than $\lambda_c \leq \bestlambdalower$. A computation by Pantone \cite{pantone-pc} using the method of differential approximants (see, e.g., \cite{GJ-new}) on the sequence $(c_n)_{n=1}^{60}$ suggests $\mu_{\rm taxi} \in [1.57376, 1.57378]$, so the limit of the present approach may in fact be $5.134$. Note that a similar situation exists for lower bounds on $\lambda_c$: 
Sinclair, Srivastava, \v{S}tefankovi\v{c} and Yin \cite{sssy} showed $\lambda_c\ge 2.538$, but Vera, Vigoda and Yang \cite{VeraVigodaYang} observed that the methods used in \cite{sssy,VeraVigodaYang} are unlikely to prove anything better than $\lambda_c\ge 3.4$, as strong spatial mixing is known not to hold at that point.

\section{Proof of phase coexistence (Theorem \ref{thm-phasecoexistence-taxi})} \label{sec-proof_phase_coexistence}

Let $\lambda > \mu_{\rm taxi}^4-1$ be fixed.  Our argument will depend on a parameter $m=m(\lambda)$ whose value will be specified later. 

Let ${\mathcal E}$ denote the set of {\em even} vertices of $\Z^2$ --- those vertices $(x,y)$ with $x+y$ even --- and let ${\mathcal O}$ denote the complementary set of odd vertices; note that these are both independent sets. 
Let $U_n$ be the box $\{-n,-(n-1), \ldots, n-1,n\}^2$. 
When $I$, $J$ are independent sets such that $I \subseteq J$, we say that $J$ {\em extends} $I$.
Let ${\mathcal J}_n^{\rm e}$ be the set of independent sets that extend ${\mathcal E}\setminus U_n$, and let $\mu_n^{\rm e}$ be the probability distribution supported on ${\mathcal J}_n^{\rm e}$ in which each set $I$ is selected with probability proportional to $\lambda^{|I \cap U_n|}$. Define $\mu_n^{\rm o}$ analogously (with ``even'' everywhere replaced by ``odd'').

Say that an independent set $I$ in ${\mathbb Z}^2$ is {\it $m$-even} if for every even $x \in U_m$, none of the four neighbors of $x$ is in $I$, and define {\it $m$-odd} analogously. Say that $I$ is {\it $m$-homogeneous} if it is either $m$-odd or $m$-even.

Now fix $n > m$ and let
${\mathcal E}_m$ be the event that an independent set is $m$-even, ${\mathcal O}_m$ the event that it is $m$-odd, and ${\mathcal H}_m$ the event that it is $m$-homogeneous. Note that all of these events are in the cylinder $\sigma$-algebra. We will establish the following conditional probability inequality for all $n > m$ and $m$ large enough:
\begin{equation} \label{distinguishing-statistic-0}
\mu_n^{\rm e}({\mathcal O}_m|{\mathcal H}_m) < 1/3.
\end{equation}
Reversing the roles of odd and even throughout the proof, we will also get
$$
\mu_n^{\rm o}({\mathcal E}_m|{\mathcal H}_m) < 1/3
$$
and so
$$
\mu_n^{\rm o}({\mathcal O}_m|{\mathcal H}_m) > 2/3.
$$
It follows that if $\mu^{\rm e}$ is any Gibbs measure obtained as a weak subsequential limit of the $\mu_n^{\rm e}$'s, and $\mu^{\rm o}$ is any obtained from the $\mu_n^{\rm o}$'s, then
$$
\mu^{\rm e}({\mathcal O}_m|{\mathcal H}_m) \leq 1/3 < 2/3 \leq \mu^{\rm o}({\mathcal O}_m|{\mathcal H}_m);
$$   
consequently, $\mu^{\rm e}$ and $\mu^{\rm o}$ are distinct Gibbs measures.

Write ${\mathcal B}_n^{\rm e}$ for the set of independent sets in $\Z^2$ that extend ${\mathcal E}\setminus U_n$ and are $m$-odd, and write ${\mathcal A}_n^{\rm e}$ for the set of independent sets in $\Z^2$ that extend ${\mathcal E}\setminus U_n$ and are $m$-homogeneous (so ${\mathcal B}_n^{\rm e} \subseteq {\mathcal A}_n^{\rm e} \subseteq {\mathcal J}_n^{\rm e}$). 

For $I \in {\mathcal A}_n^{\rm e}$ set $w_\lambda(I) = \lambda^{|I \cap U_n|}$, and for a set ${\mathcal C}$ of independent sets in ${\mathcal A}_n^{\rm e}$ let $w_\lambda({\mathcal C})$ denote $\sum_{I \in {\mathcal C}} w_\lambda(I)$. To establish (\ref{distinguishing-statistic-0}) it is enough to show
\begin{equation} \label{eq-flow_inequality}
\frac{\mu_n^{\rm e}({\mathcal O}_m)}{\mu_n^{\rm e}({\mathcal H}_m) } = \frac{w_\lambda({\mathcal B}_n^{\rm e})}{w_\lambda({\mathcal A}_n^{\rm e})} < 1/3.
\end{equation}
The intuition here is that if $I$ is conditioned to agree with ${\mathcal E}$ outside $U_n$, then under the extra condition that it is $m$-homogeneous it is far less likely to be $m$-odd than $m$-even.

We will establish (\ref{eq-flow_inequality}) by constructing, for each $I \in  {\mathcal B}_n^{\rm e}$, a collection $\varphi(I) \subseteq {\mathcal A}_n^{\rm e}$, together with a flow function $f(I,J)$ supported on $\{(I,J):I \in {\mathcal B}_n^{\rm e}, J \in \varphi(I)\}$ that satisfies
\begin{equation} \label{eq-flow_1}
\sum_{J \in \varphi(I)} f(I,J) = 1 ~~~\mbox{for each $I \in {\mathcal B}_n^{\rm e}$}
\end{equation} 
and 
\begin{equation} \label{eq-flow_2}
\sum_{I:J \in \varphi(I)} \lambda^{|I\cap U_n|-|J\cap U_n|}f(I,J) < 1/3 ~~~\mbox{for each $J \in {\mathcal A}_n^{\rm e}$}.
\end{equation}
This gives the inequality in (\ref{eq-flow_inequality}) via
\begin{eqnarray*}
w_\lambda({\mathcal B}_n^{\rm e}) & = & \sum_{I \in {\mathcal B}_n^{\rm e}} \lambda^{|I\cap U_n|} \\
& = & \sum_{I \in {\mathcal B}_n^{\rm e}} \sum_{J \in \varphi(I)} \lambda^{|I\cap U_n|}f(I,J) \\
& = & \sum_{J \in {\mathcal A}_n^{\rm e}} \lambda^{|J\cap U_n|} \sum_{I:J \in \varphi(I)} \lambda^{|I\cap U_n|-|J\cap U_n|}f(I,J) \\
& < & w_\lambda({\mathcal A}_n^{\rm e})/3. 
\end{eqnarray*}

To construct $\varphi$ we will use the fact that $I \in {\mathcal B}_n^{\rm e}$ is in {\em even phase} (predominantly even-occupied) outside $U_n$, but because $I$ is $m$-odd, it is not in even phase close to $U_m$; so there must be a {\em contour} --- an unoccupied ring of vertices --- marking the extent of the even phase inside $U_n$.

We will proceed in two stages. In Section \ref{subsec-contours}, we explain how such a contour can be explicitly constructed, and establish the various properties of the construction that we will need. (To aid readability we defer proofs of many of these properties to Sections \ref{sec:proofs-cont-facts} and \ref{sec:contours-taxi}.) 
In Section \ref{subsec-Peierls} we describe and analyze the standard Peierls argument, which involves modifying $I$ inside the contour to create $\varphi(I)$ satisfying (\ref{eq-flow_1}) and (\ref{eq-flow_2}) (for suitable choice of the flow function $f$), showing that being $m$-odd is unlikely, conditioned on being $m$-homogeneous, under even boundary condition.

\subsection{The contour and its properties} \label{subsec-contours}

Fix $I \in {\mathcal B}_n^{\rm e}$ (so $I$ includes all even vertices outside $U_n$, and is $m$-odd). Let $I'$ consist of $I$ together with each odd vertex that has none of its neighbors in $I$; notice that $I' \in {\mathcal B}_n^{\rm e}$, that it is  completely determined by $I$, and that it includes all odd vertices of $U_m$ (since $I$ is $m$-odd). The point of passing from $I$ to $I'$ is that in doing so, we ensure that the contour we construct fully encircles $U_m$ and hence has length at least on the order of $m$. We now describe how to associate with $I'$ a set $\gamma(I)$ of edges of ${\mathbb Z}^2$, which we will refer to as the {\em contour} associated with $I$. The same construction was used in \cite{bc7} and \cite{galvin}, and a very similar construction appeared in \cite{gk}.

We begin with a brief reminder of some graph theory notation. Given $S \subseteq {\mathbb Z}^2$ the subgraph of ${\mathbb Z}^2$ induced by $S$ is the graph with vertex set $S$ in which two vertices are adjacent if and only if there are adjacent in ${\mathbb Z}^2$. We will abuse notation somewhat and refer to this graph simply as $S$. Given distinct vertices $u, v \in S$ a {\em path} in $S$ from $u$ to $v$ (a $u$-$v$ path) is a sequence $u=u_0,u_1,\ldots, u_k=v$ of distinct vertices of $S$ with $u_i$ and $u_{i+1}$ adjacent for each $i=0,\ldots, k-1$. Given a vertex $u \in S$ and a subset $T \subseteq S$ of vertices with $u \not \in T$ a path in $S$ from $u$ to $T$ (a $u$-$T$ path) is any path from $u$ to $v$ with $v \in T$. A {\em component} of $S$ is a subset $C \subseteq S$ of vertices with the property that for any pair of distinct vertices $u, v \in C$ there is a path in $S$ from $u$ to $v$, and that is maximal with respect to this property.

Let $(I^{\mathcal O})^+$ be the set of odd vertices in $I'$ together with their neighbors. By the $m$-oddness of $I$, $U_m$ is contained in a single component of the graph induced by $(I^{\mathcal O})^+$; let $R$ be that component. Note that because $I'$ extends ${\mathcal E} \setminus U_n$, 
$R$ is finite, and specifically $R \subseteq U_n$.

The following property, which says that in leaving $R$ one always goes from an unoccupied even vertex to an unoccupied odd vertex (in fact, an odd vertex outside $I'$), is evident from the construction of $R$.
\begin{equation} \label{parity}
\mbox{If $uv$ is an edge with $u \in R$ and $v \not \in R$ then $u \in {\mathcal E} \setminus I$ and $v \in {\mathcal O} \setminus I'$}.
\end{equation}
Define $\gamma = \gamma(I)$ to be the set of edges $uv$ with $u \in R$ and $v \not \in R$, and such that there is a path in ${\mathbb Z}^2$ from $v$ to ${\mathbb Z}^2   \setminus   U_n$ that avoids $R$. (This last condition has the effect of removing any holes that $R$ may have.)

We now introduce a graph $\Gd$ that may be thought of as dual to $\Z^2$. The vertices of $\Gd$ are the midpoints of edges in $\Z^2$, and two such vertices are adjacent if the associated edges in $\Z^2$ are incident and perpendicular. Note that $\Gd$ is a rotated, dilated, translated copy of $\Z^2$.

Each edge in $\gamma$ is a vertex in $\Gd$, and so we may specify $\gamma$ by specifying a subgraph of $\Gd$ whose vertex set is $\gamma$. In what follows we describe a way to specify one such subgraph, which will turn out to be a cycle (a connected $2$-regular graph) in $\Gd$
that separates $R$ from ${\mathbb Z}^2 \setminus U_n$,
and moreover has a particular structure related to taxi walks.

To construct the subgraph, which we will call $\Gamma$ and also sometimes refer to as the contour associated with $I$, consider an arbitrary $uv \in \gamma$ with $u \in R$ and $v \not \in R$ (so $u \in {\mathcal E}$ and $v \in {\mathcal O}$). The edge $uv$ forms a side of two $1$-by-$1$ squares in ${\mathbb Z}^2$. Let $uvst$ be one such, with $s \in {\mathcal E}$ and $t \in {\mathcal O}$ (so each of $uv$, $vs$, $st$ and $tu$ are edges of ${\mathbb Z}^2$). Viewed as vertices of $\Gd$, $uv$ has two neighbors among the edges $uv$, $vs$, $st$ and $tu$, namely $tu$ and $vs$. 
If $s \in R$, then the directed edge $uv \rightarrow vs$ is added to $\Gamma$, and if $s \not\in R$ then the directed edge $uv \rightarrow tu$ is added to $\Gamma$. We proceed in the same way with the other $1$-by-$1$ square in ${\mathbb Z}^2$ that $uv$ forms one side of. 
(Thus $\Gamma$ is initially a directed graph; presently we will modify it slightly to create an undirected graph.)

It will be helpful for subsequent arguments to view the construction above via the following case-by-case analysis, which is easily seen to be equivalent.
There are four cases, depending on the statuses of $s$ and $t$ with regards to membership of $R$.

{\bf Case i}, $t \in R$ and $s \not \in R$: This case cannot occur since $t$ is odd and so if $t \in R$ then also~$s \in R$.

{\bf Case ii}, $t \not \in R$ and $s \not \in R$: In this case, of $tu$ and $vs$ only $tu$ is in $\gamma$ (to see that $tu \in \gamma$ note that there is a path in ${\mathbb Z}^2$ from $t$ to ${\mathbb Z}^2   \setminus   U_n$ that avoids $R$, that starts $tsv$ and then continues along any $R$-avoiding path from $v$ to ${\mathbb Z}^2   \setminus   U_n$), and we put the directed edge (in $\Gd$) $uv \rightarrow tu$ in $\Gamma$.

{\bf Case iii}, $t \in R$ and $s \in R$: In this case, of $tu$ and $vs$ only $vs$ is in $\gamma$, and we put the edge $uv \rightarrow vs$ in $\Gamma$.

{\bf Case iv}, $t \not \in R$ and $s \in R$: In this case, $vs$ is evidently in $\gamma$, and $ut$ may or not be (depending on whether there is an $R$-avoiding path in ${\mathbb Z}^2$ from $t$ to ${\mathbb Z}^2   \setminus   U_n$), and we put the edge $uv \rightarrow vs$ in $\Gamma$ (but {\em not} the edge from $uv$ to $tu$, even if $tu \in \gamma$).

What we have constructed so far is a {\em directed} graph on the set of vertices in $\Gd$ corresponding to edges in $\gamma$.
Observe that every vertex of this directed graph has out-degree two, since for each edge $uv \in \gamma$ 
exactly two edges of $\Gd$ are included: one in each of the two $1$-by-$1$ squares in ${\mathbb Z}^2$ of which $uv$ is a side. But notice that the construction of directed edges is symmetric: if we put the edge $uv \rightarrow tu$ (say) in $\Gamma$ then the construction also mandates putting in the edge $tu \rightarrow uv$. So underlying the $2$-out-regular directed graph is a $2$-regular undirected graph (a union of cycles) and it is this we take as $\Gamma$. Note that $\Gamma$ determines $\gamma$, since its vertex set is exactly $\gamma$; but $\Gamma$ is not (necessarily) the subgraph of $\Gd$ induced by $\gamma$, because 
not all edges from $\Gd$ with endpoints in $\gamma$ are in $\Gamma$ (see Case iv above). Note also that because $\Gamma$ is $2$-regular, the size of both its vertex set and its edge set is $|\gamma|$. 

If we draw $\Gamma$ in ${\mathbb Z}^2$, using straight-line segments joining midpoints of edges of $\gamma$ to represent edges of $\Gamma$, then each component of $\Gamma$ is a simple closed ${\mathbb Z}^2$-avoiding polygon in ${\mathbb R}^2$, and so encloses a finite interior (with an infinite exterior). We refer to the vertices of ${\mathbb Z}^2$ that are in the interior as the {\em vertex interior} of the component, and to all other vertices of ${\mathbb Z}^2$ as the {\em vertex exterior}. A basic fact is the following.
\begin{lemma} \label{lem-Gamma-basic_1}
If $uv \in \gamma$ with $u \in R$ and $v \not \in R$, then $u$ is in the vertex interior of the component of $uv$ in $\Gamma$, and $v$ is in the vertex exterior.
\end{lemma}

\medskip

\noindent {\bf Proof}:  
Assume without loss of generality that $uv$ is parallel to the $y$-axis. By construction $\Gamma$ has an edge to the right of $uv$ that ends at the midpoint of $uv$, and one to the left that starts at that point, so the edge $uv$ crosses $\Gamma$. This establishes that one of $u, v$ is in the vertex interior of the component of $uv$ in $\Gamma$, and the other is in the vertex exterior. If $v$ is in the vertex interior then any $v$-${\mathbb Z}^2   \setminus   U_n$ path in ${\mathbb Z}^2$ must cross $\Gamma$ and so meet $R$, contradicting the construction of $\gamma$. So $v$ is in the vertex exterior and $u$ is in the vertex interior.
\qed

\medskip

So far we have established that $\Gamma$ is a union of cycles;
Lemma \ref{lem-Gamma-basic_1} is a key ingredient in proving our first important fact about $\Gamma$; namely, that it has a single component.
\begin{lemma} \label{lem-Gamma-is-a-cycle}
The graph $\Gamma$ is a cycle. 
\end{lemma}

\medskip

\noindent {\bf Proof}: Assume, for a contradiction, that $\Gamma$ has distinct components $C_1$ and $C_2$. It cannot be the case that one of these, $C_1$ say, encloses the other. For if $uv \in \gamma$ crosses $C_2$, with $v \not \in R$, then there must (by the definition of $\gamma$) be a $v$-${\mathbb Z}^2   \setminus   U_n$ path in ${\mathbb Z}^2$ that avoids $R$; but any such path must cross $C_1$, and so meet $R$, a contradiction.

So it must be the case that $C_1$ and $C_2$ have disjoint vertex interiors. Let $uv$ cross $C_1$ and $u'v'$ cross $C_2$, with $u$ and $u'$ the interior vertices. Via Lemma \ref{lem-Gamma-basic_1} this gives an immediate contradiction: every $u$-$u'$ path in ${\mathbb Z}^2$ must use an edge of $\gamma$ and so leave $R$, contradicting the connectivity of $R$.
\qed

\medskip

The following facts about the structure of $\Gamma$ (Lemmas \ref{lem-contours_are_long4} and \ref{lem-count_of_contours}) will be used to complete the proof of Theorem \ref{thm-phasecoexistence-taxi}. We defer the proof of Lemma \ref{lem-contours_are_long4} to Section \ref{sec:proofs-cont-facts}, and that of Lemma \ref{lem-count_of_contours}, which crucially depends on the connection between contours and taxi walks (essentially, a contour is a closed taxi walk) to Section \ref{sec:contours-taxi}.
\begin{lemma} \label{lem-contours_are_long4}
$|\Gamma| \geq 2\sqrt{2}m$ and is a multiple of $4$.
\end{lemma}

Let ${\mathcal C}_\ell^m$ be the collection of all $\Gamma$ with $|\Gamma|=4\ell$ that arise in the above-described construction, as $I$ runs over ${\mathcal B}^e_n$, and let ${\mathcal C}^m = \cup_\ell {\mathcal C}_\ell^m$. 
\begin{lemma} \label{lem-count_of_contours}
There is a function $g(\ell)$ that grows subexponentially such that $|{\mathcal C}_\ell^m| \leq g(\ell)\mu_{\rm taxi}^{4\ell}$, where $\mu_{\rm taxi}$ is the connective constant of taxi walks.
\end{lemma} 

\medskip

The Peierls argument that we will use in Section \ref{subsec-Peierls} involves the {\em shift} operation. Essentially this is a shifting, by one lattice unit, of all the vertices of $I$ that are enclosed by $\Gamma$, while leaving the remainder of $I$ unchanged. The content of Lemma \ref{lem-shift} below is that this allows $I$ to be augmented substantially, leading to a weight-increasing map from ${\mathcal B}^e_n$ to ${\mathcal A}^e_n$ that allows for the construction of a flow function satisfying (\ref{eq-flow_1}) and (\ref{eq-flow_2}).

Let $W$ be the vertex interior of $\Gamma$ and $W'$ the vertex exterior. For $v \in {\mathbb Z}^2$, and $s \in \{(1,0), (0,1), (-1,0), (0,-1)\}$, let $\sigma_s(v)=v+s$. Let $I_s = (I \cap W') \cup \{\sigma_s(v): v \in I \cap W\}$ --- $I_s$ may be thought of as the result of shifting $I$ one unit in the $s$ direction within $W$, while leaving it unchanged outside $W$. Let $\tilde{I}_s = \{v \in W: \sigma_s^{-1}(v) \in W'\}$ --- we may think of $\tilde{I}_s$ as the set of vertices in $W$ with the property that their preimage under the map $\sigma_s$ is outside $W$; note that $\sigma_s^{-1}(v)=v-s$. Finally, let $I''_s = I_s \cup \tilde{I}_s$. 
\begin{lemma} \label{lem-shift}
The shifted set $I_s$ is an independent set, with $|I_s|=|I|$. Moreover, the augmented shifted set $I''_s$ is an independent set, with $|I''_s|=|I_s|+|\tilde{I}_s|$. Finally, there is a choice of $s$ for which $|\tilde{I}_s| \geq |\gamma|/4$.
\end{lemma}

\medskip

The proof of Lemma \ref{lem-shift} uses standard ideas (see, e.g., \cite[Proposition 2.12]{gk}, \cite[Lemma 4.1]{galvin} and \cite[proof of Lemma 6]{bc7}). The same is true for our next lemma (see, e.g., \cite[equation (15)]{gk}). For completeness we furnish proofs in Section \ref{sec:proofs-cont-facts}. 
\begin{lemma} \label{lem-reconstruct}
If $I \in {\mathcal B}_n^{\rm e}$ has associated contour $\Gamma$, and $J = I_s \cup S$ where $S \subseteq \tilde{I}_s$ and $s$ is one of the four possible shift directions, then $I$ is completely determined by $J$, $s$ and $\gamma$.
\end{lemma}

\medskip

We also need one more lemma, that says that after shifting we go from $m$-odd independent set to an $m$-homogeneous set. The proof appears in Section \ref{sec:proofs-cont-facts}. 
\begin{lemma} \label{lem-good_map}
If $I \in {\mathcal B}_n^{\rm e}$ then $I''_s \in {\mathcal A}_n^{\rm e}$.
\end{lemma}

\subsection{The Peierls argument for phase coexistence} \label{subsec-Peierls}

We are now in a position to define
the collection $\varphi(I) \subseteq {\mathcal A}_n^{\rm e}$ and $f(I,J)$ for $I \in {\mathcal B}_n^{\rm e}$ and $J \in \varphi(I)$. First, choose an $s \in \{(1,0), (0,1), (-1,0), (0,-1)\}$ for which $|\tilde{I}_s| \geq |\gamma|/4$ (by Lemma \ref{lem-shift} there is such an $s$; choose, for example, the first such that works in some arbitrary ordering). Next, set 
$$
\varphi(I) = \{I_s \cup S:S \subseteq \tilde{I}_s\}
$$
(by Lemma \ref{lem-good_map} we have $\varphi(I) \subseteq {\mathcal A}_n^{\rm e}$). Finally, for $I \in {\mathcal B}_n^{\rm e}$ and $J \in \varphi(I)$ set 
$$
f(I,J) = \frac{\lambda^{|S|}}{(1+\lambda)^{|\tilde{I}_s|}}.
$$
For this choice of $f$ we have
$$
\sum_{J \in \varphi(I)} f(I,J) = \sum_{S \subseteq \tilde{I}_s} \frac{\lambda^{|S|}}{(1+\lambda)^{|\tilde{I}_s|}} = \frac{1}{(1+\lambda)^{|\tilde{I}_s|}} \sum_{k=0}^{|\tilde{I}_s|} {|\tilde{I}_s| \choose k} \lambda^k = 1,
$$
so (\ref{eq-flow_1}) is established, and it only remains to verify (\ref{eq-flow_2}). 

We now present the Peierls argument that verifies (\ref{eq-flow_2}) and completes the proof of phase coexistence. Recall that $\lambda > \mu_{\rm taxi}^4-1$ has been given. Choose $\mu>0$ to be such that $\mu^4-1$ is the midpoint of $[\mu_{\rm taxi}^4-1,\lambda]$. Choose $m$ sufficiently large that 
\begin{equation} \label{m-cond_1}
g(\ell)\mu_{\rm taxi}^{4\ell} < \mu^{4\ell}
\end{equation}
for all $\ell \geq \sqrt{2} m/2$, where $g(\ell)$ 
is the subexponential function from Lemma \ref{lem-count_of_contours} (this inequality holds for all sufficiently large $m=m(\lambda)$ since $\mu > \mu_{\rm taxi}$, and depends only on $\lambda$), and also that
\begin{equation} \label{m-cond_2}
\sum_{\ell \geq \sqrt{2} m/2} \frac{\mu^{4\ell}}{(1+\lambda)^{\ell}} < 1/3
\end{equation}
(this inequality holds for all sufficiently large $m=m(\lambda)$ since $1+\lambda > \mu^4$).

For $n > m$, fix $J \in {\mathcal A}_n^{\rm e}$. From the definitions of $\varphi(I)$ and $f$, we have that
\begin{equation}
	\sum_{I:J \in \varphi(I)} \lambda^{|I|-|J|}f(I,J) 
	 =   \sum_{I:J \in \varphi(I)}\frac{1}{(1+\lambda)^{|\tilde{I}_s|}}
	\le \sum_{I:J \in \varphi(I)}\frac{1}{(1+\lambda)^{|\gamma(I)|/4}}, \label{int1}		
\end{equation}
since $|\tilde{I}_s| \geq |\gamma(I)|/4$ by Lemma \ref{lem-shift}. Bearing Lemma \ref{lem-reconstruct} in mind, for each $\Gamma \in {\mathcal C}^m$ there are at most four $I \in {\mathcal B}_n^{\rm e}$ for which $J \in \varphi(I)$ (at most one for each shift $s \in \{(1,0), (0,1), (-1,0),$ $(0,-1)\}$). Hence, 
	\begin{eqnarray}		
	\sum_{I:J \in \varphi(I)} \lambda^{|I|-|J|}f(I,J)  & \leq & 4 \sum_{\Gamma \in {\mathcal C}^m} \frac{1}{(1+\lambda)^{|\Gamma|/4}} \nonumber \\
	& = & 4 \sum_{4\ell \geq 2\sqrt{2} m} \frac{|{\mathcal C}_\ell^m|}{(1+\lambda)^{\ell}} \label{long_contours_X} \\
	& \leq & \sum_{\ell \geq \sqrt{2} m/2} \frac{g(\ell)\mu_{\rm taxi}^{4\ell}}{(1+\lambda)^{\ell}} \nonumber \\
	& < & \sum_{\ell \geq \sqrt{2} m/2} \frac{\mu^{4\ell}}{(1+\lambda)^{\ell}} \nonumber \\
	& < & 1/3, \nonumber
	\end{eqnarray} 
where the first inequality uses (\ref{int1}), we use Lemma \ref{lem-contours_are_long4} in (\ref{long_contours_X}), the second inequality uses Lemma \ref{lem-count_of_contours}, and the remaining inequalities follow from our choice of $m$ (specifically using (\ref{m-cond_1}) and (\ref{m-cond_2})). The proof is now complete.
\qed

\subsection{Proofs of contours facts}
\label{sec:proofs-cont-facts}

In this section we begin wrapping up the proof of Theorem \ref{thm-phasecoexistence-taxi} by providing the proofs of Lemmas \ref{lem-contours_are_long4}, \ref{lem-shift}, \ref{lem-reconstruct} and \ref{lem-good_map}.

We will need a basic fact about $\Gamma$ that comes immediately from the construction.
\begin{lemma} \label{lem-Gamma-basic_2}
If $\{a,b,c,d\}$ are the vertices of a $1$-by-$1$ square in ${\mathbb Z}^2$ (with $ab$, $bc$, $cd$ and $da$ the edges of ${\mathbb Z}^2$), then in $\Gamma$ it is not possible for $bc$ to be adjacent to both $ab$ and $cd$.
\end{lemma}

\medskip

\noindent {\bf Proof of Lemma \ref{lem-contours_are_long4}}:  
Since the interior of $\Gamma$ contains a vertex with $x$-coordinate $m$ (along the top of $U_m$) and one with $x$-coordinate $-m$ (along the bottom), and each edge of $\Gamma$ spans a distance of $1/\sqrt{2}$ in the $x$-direction, it follows that $\Gamma$ must have at least $2\sqrt{2}m$ edges.

To argue about the length of $\Gamma$, we view it as a simple closed ${\mathbb Z}^2$-avoiding polygon in ${\mathbb R}^2$, in the manner described before the statement of Lemma \ref{lem-Gamma-basic_1}, and consider traversing this polygon in a clockwise direction starting from an arbitrarily chosen point $P$ that is the midpoint of the form $(a_x, a_y-1/2)$ of an edge of $\gamma$, with $a_x$, $a_y$ integers. We traverse in steps of length $\sqrt{2}/2$, which corresponds to moving from the midpoint of one edge of $\gamma$ to the midpoint of an adjacent (and perpendicular) edge. 

A complete traverse of the polygon consists of $x_\searrow$ steps oriented southeast (parallel to the edge from $(0,0)$ to $(-1,-1)$), $x_\nwarrow$ steps oriented northwest, $x_\swarrow$ steps oriented southwest and $x_\nearrow$ steps oriented northeast, and because $\Gamma$ is closed we have $x_\searrow=x_\nwarrow$ and $x_\swarrow=x_\nearrow$.

Starting at $P$, a point in ${\mathbb R}^2$ of the form $(a_x, a_y-1/2)$ with $a_x$, $a_y$ integers, after two steps we return to a point of this form, having passed through a point of the form $(a'_x-1/2, a'_y)$ with $a'_x$, $a'_y$ integers. These two steps must be one of: southwest followed by southeast or vice versa; northwest followed by northeast or vice versa; or two steps in the same direction. (All other possibilities, such as southwest followed by northwest, are ruled out by Lemma \ref{lem-Gamma-basic_2}).

Write $x_{\swarrow\searrow}$ for the total number (over the entire polygon) of pairs of steps of the kind just described that consist of southwest followed by southeast, and write $x_{\searrow\swarrow}$, $x_{\nwarrow\nearrow}$, $x_{\nearrow\nwarrow}$, $x_{\searrow\searrow}$, $x_{\nwarrow\nwarrow}$, $x_{\swarrow\swarrow}$ and $x_{\nearrow\nearrow}$ for the count of the other possible pairs. Using $x_\searrow=x_\nwarrow$ we get
$$
x_{\swarrow\searrow}+x_{\searrow\swarrow}+2x_{\searrow\searrow} = x_{\nwarrow\nearrow}+x_{\nearrow\nwarrow}+2x_{\nwarrow\nwarrow}
$$      
and using $x_\swarrow=x_\nearrow$ we get
$$
x_{\searrow\swarrow}+x_{\swarrow\searrow}+2x_{\swarrow\swarrow} = x_{\nwarrow\nearrow}+x_{\nearrow\nwarrow}+2x_{\nearrow\nearrow}.
$$
Combining (and dividing by 2) we get
$$
x_{\swarrow\searrow}+x_{\searrow\swarrow}+x_{\searrow\searrow}+x_{\swarrow\swarrow} = x_{\nwarrow\nearrow}+x_{\nearrow\nwarrow}+x_{\nwarrow\nwarrow}+x_{\nearrow\nearrow}.
$$
It follows that $x_{\swarrow\searrow}+x_{\searrow\swarrow}+x_{\searrow\searrow}+x_{\swarrow\swarrow} + x_{\nwarrow\nearrow}+x_{\nearrow\nwarrow}+x_{\nwarrow\nwarrow}+x_{\nearrow\nearrow}$ is even, and so $|\Gamma|$, being twice this sum, is a multiple of $4$.
\qed

That $4$ divides $|\Gamma|$ could also be read out of \cite[Lemma 5]{bc7}; we give a self-contained proof above to avoid a lengthy detour matching our notation to that of \cite{bc7}. 
 
\medskip

\noindent {\bf Proof of Lemma \ref{lem-shift}}: We begin with the final statement. Let $\gamma_s$ be the set of edges in $\gamma$ of the form $uv$ with $u \in W$, $v \not \in W$ and $v = u - s$. Note that $\gamma = \cup_s \gamma_s$, so there is a choice of $s$ for which $|\gamma_s|\geq |\gamma|/4$. Now the map from $\gamma_s$ to $\tilde{I}_s$ that sends $uv$ to $u$ is injective (for each $u \in \tilde{I}_s$ there is a unique $v$ such that $uv \in \gamma_s$, namely $u-s$), so $|\tilde{I}_s| \geq |\gamma_s|\geq |\gamma|/4$.

To show $|I_s|=|I|$ consider the map $\varphi$ from $I$ to $I_s$ that sends $v$ to $v$ if $v \in W'$ and sends $v$ to $v+s$ if $v \in W$. The restrictions of $\varphi$ both to $I \cap W'$ and to $I \cap W$ are bijections. Also, $\varphi(I \cap W')$, being $I \cap W'$, is disjoint from $W$, and $\varphi(I \cap W) \subseteq W$, this latter since the vertices of $W$ with a neighbor outside $W$ are all unoccupied. This shows $|I_s|=|I|$. To see that $I_s$ is an independent set, note first that $\varphi(I \setminus W)$ and $\varphi(I \cap W)$ are both independent sets, so we need only rule out the possibility of having $v_1 \in I \cap W$ and $v_2 \in I  \cap W'$ with $\varphi(v_1)\varphi(v_2)$ an edge in ${\mathbb Z}^2$. Since (as we have already observed) $\varphi(v_1) \in W$, and $\varphi(v_2) \in W'$, such an adjacency would put $\varphi(v_1)\varphi(v_2)$ in $\gamma$; but since $\varphi(v_2)=v_2$ this would lead to an edge in $\gamma$ with one endvertex occupied, contradicting (\ref{parity}).

Next we show that $\tilde{I}_s$ is disjoint from $I_s$. It is clearly disjoint from $\varphi(I \cap W')$. It is also easily seen to be disjoint from $\varphi(I \cap W)$, since all vertices $v$ in $\varphi(I\cap W)$ have $v-s \in W$, and no vertices in $\tilde{I}_s$ have this property. Finally we need to show that no $v \in \tilde{I}_s$ is adjacent to something in $I'_s$. There cannot be a $w \in \varphi(I \cap W')$ with $vw \in {\mathbb Z}^2$, for $vw$ would then be in $\gamma$ and have one endvertex ($w$) occupied. Next we consider a $w \in \varphi(I \cap W)$ with $vw$ an edge of ${\mathbb Z}^2$. We cannot have $w=v+s$, for then we would have $v \in I$ (again creating an edge in $\gamma$ with one endvertex, this time $v$, occupied). We cannot have $w=v-s$ since this would put $w$ into $W'$ (by definition of $\tilde{I'}_s$), and we know $w \in W$ since $\varphi(I' \cap W) \subseteq W$. There remains the case $w=v+s'$, with $s'$ perpendicular to $s$. But in this case, $w-s$ is an occupied vertex in $W$, and $v-s$ is a neighbor of $w-s$ that is outside $W$, again creating an impossible edge in $\gamma$. \qed

\medskip

\noindent {\bf Proof of Lemma \ref{lem-reconstruct}}: $\Gamma$ determines $I$ and so $W$, and this together with $s$ determines $\tilde{I}_s$ (which, crucially, depends only on $W$ and $s$ and not on $I$). This allows $S$ to be determined, as $S=J \cap \tilde{I}_s$, from which $I_s$ can be determined as $I_s=J\setminus S$. Finally we determine $I$ as $I= (I_s \cap W') \cup \{\sigma_{s}^{-1}(v): v \in I_s \cap W\}$.
\qed

\medskip

\noindent {\bf Proof of Lemma \ref{lem-good_map}}: 
Because $I$ is $m$-odd we know that no even vertex of $U_m$ is in $I$, and nor is any even vertex outside $U_m$ that is adjacent to something in $U_m$. We aim to establish that after the shift operation no odd vertex of $U_m$ is in $I$, and nor is any odd vertex outside $U_m$ that is adjacent to something in $U_m$; this shows that $I_s$ is $m$-even, and since in going from $I_s$ to $I''_s$ we only add even vertices, so also is $I''_s$.

That no odd vertex of $U_m$ is in $I$ after the shift is clear, since $R$ includes $\sigma_{-s}(v)$ for every odd $v \in U_m$, no such $\sigma_{-s}(v)$ is in $I$, and the status of $v$ with regards membership of $I_s$ is identical to the status of $\sigma_{-s}(v)$ with regards membership of $I$. The same argument holds for any odd vertex $v$ outside $U_m$ adjacent to something in $U_m$ for which either $\sigma_{-s}(v) \in U_m$ or  $\sigma_{-s}(v)$ is outside $U_m$ but adjacent to something in $U_m$. 

There remains the case of odd $v$, adjacent to something in $U_m$, with $\sigma_{-s}(v)$ not in $U_m$ and not adjacent to something in $U_m$. If such a $v$ is not in $I$, then it is clearly not in $I_s$. If $v \in I$ then $v \in R$ and $\sigma_{-s}(v) \not \in I$ and so as before $v \not \in I_s$.
\qed

\subsection{Contours as taxi walks} \label{sec:contours-taxi}

In this section we establish the connection between contours and taxi walks, which allows us to give the proof of Lemma \ref{lem-count_of_contours}. The key ingredient is the following. 
\begin{lemma} \label{lem-Gamma-is-a-taxi-walk}\ 
\begin{enumerate}
\item Viewed as a polygon in $\Gd$, if $\Gamma$ turns, goes straight for an odd number of steps, and turns again, then the second turn must be in the same direction as the first, while if it goes straight for an even number of steps, then the second turn must be in the opposite direction.
\item If $\Gamma$ turns, then it cannot turn again after a single step.
\end{enumerate}
\end{lemma}

\medskip

\noindent {\bf Proof}: We first show that if $\Gamma$ turns, goes straight for an odd number of steps, and turns again, then the second turn must be in the same direction as the first, while if it goes straight for an even number of steps, then the second turn must be in the opposite direction.

Suppose that $(e_1, f_1, f_2, \ldots, f_{2k+1}, e_2)$ is a list of consecutive edges in $\Gamma$, with $e_1$ perpendicular to $f_1$, all the $f_i$'s parallel, and $e_2$ perpendicular to $f_{2k+1}$.  Without loss of generality $e_1$ is the edge in $\Gamma$ from $(-1/2,1)$ to $(0,1/2)$, which by Lemma \ref{lem-Gamma-basic_2} forces $f_1$ to go from $(0,1/2)$ to $(1/2,1)$, which forces $f_i$ to go from $(0,1/2) + (i-1)(1/2,1/2)$ to $(0,1/2) + i(1/2,1/2)$), and in particular $f_{2k+1}$ to go from $(k,k+1/2)$ to $(k+1/2,k+1)$. Again by Lemma \ref{lem-Gamma-basic_2}, $e_2$ must now go from $(k+1/2,k+1)$ to $(k,k+3/2)$. This shows that two turns in $\Gamma$ separated by an odd number of steps must both go in the same direction (counterclockwise in this case). The case of $\Gamma$ taking an even number of steps between turns is dealt with similarly.

Next we consider the possibility of $\Gamma$ taking two consecutive turns. Suppose that the turns are taken around vertex $v$, in the sense that $v$ has neighbors (read off in cyclic order) $a$, $b$, $c$, $d$, and $\Gamma$ has edges from $av$ to $bv$, from $bv$ to $cv$, and from $cv$ to $dv$. (Bearing Lemma \ref{lem-Gamma-basic_2} in mind, no other situation is possible). 

Consider first the case where $v$ is an odd vertex. In the $1$-by-$1$ square that $a$, $v$ and $d$ are three corners of, the construction of $\Gamma$ dictates that there must be an edge of $\Gamma$ from $av$ to $dv$ (we are either in case iii or case iv). This means that $\Gamma$ encloses just the odd vertex $v$, which, by Lemma \ref{lem-Gamma-basic_1}, cannot happen.

If $v$ is even, let $q$ be the vertex that completes the $1$-by-$1$ square that includes $b$ and $c$ as corners, $r$ the one for $c$ and $d$, $s$ the one for $a$ and $b$, and $t$ the one for $a$ and $d$. Note that $\Gamma$ cannot have an edge from $av$ to $dv$, for if it did it would be a $4$-cycle enclosing a single vertex $v$, implying that $|R|=1$, a contradiction since $U_m \subseteq R$. Looking at the construction rules for $\Gamma$, we see that we must have $q \not \in R$, $r \not \in R$, $s \not \in R$ and $t \in R$, and $\Gamma$ has edges from $av$ to $at$ and from $dv$ to $dt$. Note now that $v \in R$ has all four of its neighbors outside $R$. This is not possible, since by construction of $R$ every even vertex of $R$ must have a neighbor in $R$. It follows that $\Gamma$ cannot take two consecutive turns. \qed

\medskip

The import of Lemma \ref{lem-Gamma-is-a-taxi-walk} is that $\Gamma$ can be thought of as a taxi walk; we now record this key fact formally. 
\begin{lemma} \label{lem-Gamma-is-a-taxi-walk_2}
Let $m_x, m_y$ be integers such that $(m_x, m_y-1/2)$ is the apex of a ``vee'' in $\Gamma$; that is to say, $(m_x-1/2,m_y)$ is adjacent to $(m_x, m_y-1/2)$ in $\Gamma$, and $(m_x, m_y-1/2)$ is adjacent to $(m_x+1/2,m_y)$. There is a unique orientation of the edges of $\Gd$ such that it becomes isomorphic to $\Ztaxi$ via a translation that sends $(m_x, m_y-1/2)$ to the origin, followed by a clockwise rotation through $\pi/4$, followed by a dilation by $\sqrt{2}$. Under this orientation, if the edge from $(m_x-1/2,m_y)$ to $(m_x, m_y-1/2)$ is removed from $\Gamma$ then the residue is mapped to a taxi walk of length $|\Gamma|-1$.    
\end{lemma}

\medskip

\noindent {\bf Proof of Lemma \ref{lem-count_of_contours}}: By Lemma \ref{lem-Gamma-is-a-taxi-walk_2} an element $\Gamma$ of ${\mathcal C}_\ell^m$ is fully described by specifying a midpoint $(m_x, m_y-1/2)$ ($m_x, m_y$ integers) of an edge in ${\mathbb Z}^2$ where $\Gamma$ makes a ``vee'' turn, followed by specifying a taxi walk of length $4\ell-1$. Since $\Gamma$ is a simple closed curve of length $4\ell/\sqrt{2}$ that encloses the origin, there are at most $O(\ell^2)$ choices for the pair $(m_x,m_y)$, and by Lemma \ref{thm-best_mut} there are at most $f_{\rm taxi}(4\ell-1)\mu_{\rm taxi}^{4\ell-1}$ choices for the taxi walk. The lemma follows.    
\qed

\section{The taxi walk connective constant (Theorem \ref{thm-best_mut})} \label{sec-taxibounds}

In this section we prove Theorem \ref{thm-best_mut}. A helpful initial observation is that $\Ztaxi$ is vertex-transitive; specifically, for each $(x,y) \in \Z^2$, the bijective map $f_{(x,y)}:\Z^2 \rightarrow \Z^2$ given by
$$
f_{(x,y)} = \left\{
\begin{array}{ll}
\mbox{translation by $(-x,-y)$} & \mbox{if $x,y$ both even} \\
\mbox{translation by $(-x,-y)$, then rotation through $\pi$ radian} & \mbox{if $x,y$ both odd} \\
\mbox{translation by $(-x,-y)$, then reflection across $x$-axis} & \mbox{if $x$ odd, $y$ even} \\
\mbox{translation by $(-x,-y)$, then reflection across $y$-axis} & \mbox{if $x$ even, $y$ odd}
\end{array}
\right.
$$  
induces an orientation-preserving bijection of $\Ztaxi$ that sends $(x,y)$ to the origin. 

We begin the proof of Theorem \ref{thm-best_mut} by establishing the submultiplicativity of $\c_n$ (or, equivalently,  the subadditivity of $\log \c_n$).
\begin{lemma}\label{submult}  For $n, m \geq 1$, $\c_{n+m} \leq \c_n \c_m$.
\end{lemma}

\noindent {\bf Proof}: If we split a taxi walk of length $n+m$  into two pieces, an initial segment of length $n$ and a terminal segment of length $m$, then both resulting pieces are self-avoiding. Moreover the initial segment of length $n$ is a taxi walk of length $n$, while the terminal segment of length $m$ gets mapped to a taxi walk of length $m$ by the map $f_{(x,y)}$ described above, where $(x,y)$ is the initial vertex of the terminal segment. It is straightforward to verify that this gives rise to an injective mapping from taxi walks of length $n+m$ to ordered pairs of taxi walks, the first of length $n$ and the second of length $m$, so that $\c_{n+m} \leq \c_n  \c_m$.
\qed

\medskip

It follows from Lemma~\ref{submult} that $d_n := \log \c_n$ is subadditive, i.e., $d_{n+m} \leq d_n + d_m$.
By Fekete's Lemma (see, e.g., \cite[Lemma 1.2.2]{steele}) we know that $\lim_{n\rightarrow \infty} d_n /n$ exists and that
\begin{equation}\label{fekete}
\lim_{n\rightarrow \infty} \frac{d_n}{n} = \inf_n \frac{d_n}{n}.
\end{equation}
Thus we can write the number of taxi walks of length $n$ as $\c_n = f_{\rm taxi}(n) \mu_{\rm taxi}^n$, where  $\mu_{\rm taxi}$ is a constant and $f_{\rm taxi}(n)$ is subexponential in $n$.

We have already (in Section \ref{sec-taxi-walks}) observed that $\sqrt{2} \leq \mu_{\rm taxi} \leq (1+\sqrt{5})/2$. Various techniques from the self-avoiding walk literature --- subadditivity, Alm's method, the Goulden--Jackson cluster method, and Kesten's methods of bridges and irreducible bridges --- can be used to improve both bounds. We now  discuss these methods and our associated results. 

\subsection{Upper bounds on $\mu_{\rm taxi}$}

Subadditivity gives us a strategy for getting a better upper bound on $\mu_{\rm taxi}$.  From (\ref{fekete}) we see that for all $n$, $\log \c_n /n$ is an upper bound for $\log \mu_{\rm taxi}$. Then,
using that $c_{60}=2189670407434$ (see Table \ref{table:c} and \cite{BGRT-data}) gives the bound $\mu_{\rm taxi} < 1.60574$ and $\mu_{\rm taxi}^4-1 < 5.6482$. 

\begin{table}
	\begin{center}
		\begin{tabular}{|l l|l l|l l|l l|l l|}
			\hline
			$c_{1}$ & 2&$c_{13}$ & 740&$c_{25}$ & 208506&$c_{37}$ & 54807754&$c_{49}$ & 13922238632\\
			$c_{2}$ & 4&$c_{14}$ & 1192&$c_{26}$ & 332616&$c_{38}$ & 87077354&$c_{50}$ & 22069957494\\
			$c_{3}$ & 6&$c_{15}$ & 1918&$c_{27}$ & 530588&$c_{39}$ & 138346766&$c_{51}$ & 34986181158\\
			$c_{4}$ & 10&$c_{16}$ & 3064&$c_{28}$ & 843222&$c_{40}$ & 219324398&$c_{52}$ & 55383388278\\
			$c_{5}$ & 16&$c_{17}$ & 4910&$c_{29}$ & 1342662&$c_{41}$ & 348109128&$c_{53}$ & 87740467384\\
			$c_{6}$ & 26&$c_{18}$ & 7872&$c_{30}$ & 2138280&$c_{42}$ & 552582790&$c_{54}$ & 139014623272\\
			$c_{7}$ & 42&$c_{19}$ & 12620&$c_{31}$ & 3405346&$c_{43}$ & 877163942&$c_{55}$ & 220254102104\\
			$c_{8}$ & 68&$c_{20}$ & 20114&$c_{32}$ & 5406522&$c_{44}$ & 1389806294&$c_{56}$ & 348536652664\\
			$c_{9}$ & 110&$c_{21}$ & 32150&$c_{33}$ & 8597632&$c_{45}$ & 2204289314&$c_{57}$ & 551914140382\\
			$c_{10}$ & 178&$c_{22}$ & 51396&$c_{34}$ & 13674278&$c_{46}$ & 3496483316&$c_{58}$ & 874039817792\\
			$c_{11}$ & 288&$c_{23}$ & 82160&$c_{35}$ & 21748530&$c_{47}$ & 5546212122&$c_{59}$ & 1384184997874\\
			$c_{12}$ & 460&$c_{24}$ & 130730&$c_{36}$ & 34501460&$c_{48}$ & 8783360626&$c_{60}$ & 2189670407434\\
			\hline
		\end{tabular}
	\end{center}
	\caption{Values of $c_i$ for $i=1,\dots,60$; see \cite{BGRT-data}.}
	\label{table:c}
\end{table}

The connective constant for ordinary self-avoiding walks has been well studied, and some of the methods used to obtain bounds there can be adapted to deal with taxi walks. In this section, we adapt two methods due to Alm \cite{Alm} and Goulden and Jackson \cite{GJ}
to bound $\mu_{\rm taxi}$ and thus establish Theorem \ref{thm-best_mut}. 
The bounds derived using these methods are 
very similar, so both are provided for completeness. 

First we discuss the method of Alm \cite{Alm}. Fix $n > m >0$. Construct a square matrix $A(m,n)$ whose $ij$ entry counts the number of taxi walks of length $n$ that begin with the $i$th taxi walk of length $m$, and end with the $j$th taxi walk of length $m$, for some fixed ordering of the walks of length $m$ (formally we mean that if $(x,y)$ is the vertex that begins the terminal segment of length $m$ of the walk, then the map $f_{(x,y)}$ described earlier sends the the terminal segment of length $m$ to the $j$th taxi walk of length $m$). Then a result of Alm \cite[Theorem 1]{Alm} says that 
$$
\mu_{\rm taxi} \leq \lambda_1(A(m,n))^{1/(n-m)},
$$
where $\lambda_1$ indicates the largest positive eigenvalue. (Note that when $m=0$ this recovers the subadditivity bound discussed earlier). Alm's result as stated in \cite{Alm} only applies to bound the ordinary connective constant of a finitely generated lattice, directed or otherwise. His proof is easily seen to go through without change, however, when the extra condition is added that walks do not take two consecutive turns.  
We calculated $A(20,60)$; this is a square matrix of dimension 20114
for which we can estimate its largest eigenvalue using MATLAB. This gives that $\mu_{\rm taxi} < 1.58834$ and $\mu_{\rm taxi}^4-1 < 5.3646$ 
(again see \cite{BGRT-data} for this data).

Our second approach is the Goulden--Jackson cluster method \cite{GJ}. This is an algorithm which takes as input a finite alphabet $A$, an integer $n$ and a finite list ${\mathcal M}$ of words over $A$ --- the elements of which we refer to as {\em mistakes} --- and outputs the number $\ell_n$ of words of length $n$ over $A$ that do not contain any mistakes as subwords (that is, as strings of consecutive letters in the word). 

Recall from the proof of Lemma \ref{lem-fib} that a taxi walk of length $n$ may be encoded by a pair $(a,\sigma)$, where $a \in \{N,E\}$ and $\sigma$ is a word of length $n-1$ over alphabet $\{s,t\}$.  Suppose that ${\mathcal M}$ is a finite set of subwords that is not allowed to occur in any word $\sigma$ over alphabet $\{s,t\}$ that occurs in an encoding of a taxi walk (for example, $tt$ is one such subword). If $\ell_n$ is as defined in the last paragraph then we have $\c_{n+1} \leq 2\ell_n$, so by subadditivity $\mu_{\rm taxi} \leq (2\ell_n)^{1/(n+1)}$.    

We can improve this slightly. Alm \cite[Remark 9]{Alm} observes that in a vertex-transitive lattice for which any self-avoiding walk of length $1$ can be mapped on to any other self-avoiding walk of length $1$ by some orientation-preserving symmetry (built from translations, rotations and reflections), the connective constant is bounded above by $(f(n+1)/f(1))^{1/n}$, where $f(m)$ is the number of self-avoiding walks of length $m$ starting for some fixed vertex (by vertex-transitivity, it does not matter which). Applying this to the present situation (where Alm's condition is certainly satisfied, with $f(1)=2$), we get 
\begin{equation} \label{eq-GJ_Alm}
\mu_{\rm taxi} \leq \ell_n^{1/n}.
\end{equation}
 
If $(a,\sigma')$ encodes a walk in the Manhattan lattice that takes no two consecutive turns, starts and ends at the origin, and otherwise does not visit any vertex twice, then it is evident that the word $\sigma'$ cannot occur as a subword of $\sigma$ in any taxi walk $(a, \sigma)$. We refer to such a $\sigma'$ as a {\em taxi polygon} of length $|\sigma'|+1$, where $|\sigma'|$ is the number of letters in $\sigma'$. For example, $sstsstsstss$ is a taxi polygon of length 12, and $tstsstsssstsssstsst$ is a taxi polygon of length 20.

We have enumerated taxi polygons of length at most 48 (there are 8,009,144 of them). We then used an implementation of the Goulden--Jackson cluster method due to Noonan and Zeilberger \cite{NZ} to calculate $a_{802}$ with ${\mathcal M}$ consisting of the set of all taxi polygons of length at most 44 (there are 1,721,326 of them) together with the word $tt$. Via (\ref{eq-GJ_Alm}) this leads to $\mu_{\rm taxi} < \bestmuupper$ and $\mu_{\rm taxi}^4-1 < \besttorus$, as stated in Theorem \ref{thm-best_mut}. (See \cite{BGRT-data} for the data and the computer code used to generate it.) 

\subsection{Lower bounds on $\mu_{\rm taxi}$}

To improve the trivial lower bound $\sqrt{2} \leq \mu_{\rm taxi}$ we consider {\em bridges} (introduced for ordinary self-avoiding walks by Kesten \cite{Kesten}). A {\em bridge}, for our purposes, is a taxi walk that begins by moving from the origin $(0,0)$ to the vertex $(1,0)$, never revisits the $y$-axis, and ends by taking a step parallel to the $x$-axis to a vertex on the walk that has maximum $x$-coordinate over all vertices in the walk (but note that this maximum does not have to be uniquely achieved at the final vertex).

Let $b_n$ be the number of bridges of length $n$ (by convention $b_0=1$).  Observe that bridges are supermultiplicative, that is, $b_{n+m} \geq b_n  b_m$ (and $\log b_n$ is superadditive).  To see this, consider bridges $\beta_1$ of length $n$ and $\beta_2$ of length $m$. By the definition of a bridge, it is straightforward to verify that if we concatenate $\beta_1$ and the image of $\beta_2$ under the map $f^{-1}_{(x,y)}$, where $(x,y)$ is the terminal vertex of $\beta_1$, then the result is a bridge. Moreover, the map just described from pairs of bridges, the first of length $n$ and the second of length $m$, to bridges of length $n+m$, is injective.
It follows that there are at least $b_n^k$ taxi walks of length $kn$ (just concatenate $k$ length $n$ bridges), so that
$$
\mu_{\rm taxi} = \lim_{m \rightarrow \infty} \c_m^{1/m} \geq \lim_{k \rightarrow \infty} (b_n^k)^{1/nk} = b_n^{1/n}.
$$
Since $b_{60}=80312795498$ (see \cite{BGRT-data}), we get that $\mu_{\rm taxi} > 1.51965$ and $\mu_{\rm taxi}^4-1 > 4.3330$.

Using Kesten's more sophisticated notion of {\em irreducible bridges} (bridges that are not the concatenation of shorter bridges), we can get significantly better bounds. Our discussion follows the approach of Alm and Parviainen \cite{AlmParviainen}). 

Say that an internal vertex $(x,y)$ along a bridge is a {\em cutvertex} if the walk from the origin up to $(x,y)$ is a bridge, after $(x,y)$ the next vertex of the walk is $(x+1,y)$, and the walk from $(x,y)$ to the end (more correctly, the image of this walk under $f_{(x,y)}$) is also a bridge. Say that a bridge is {\em irreducible} if it does not have a cutvertex. Denote by $a_n$ the number of irreducible bridges of length $n$ (by convention $a_0=0$).

Fix $n \geq 1$. For each $\ell \geq 1$, each solution to $k_1+\ldots +k_\ell = n$ with each $k_i \geq 1$, and each sequence of irreducible bridges $(p_1, \ldots, p_\ell)$ with $p_i$ of length $k_i$ for each $i$, there corresponds a bridge of length $n$ obtained by concatenating the $p_i$'s (with suitable translations, reflections and rotations where necessary). Moreover, each bridge of length $n$ is obtained exactly once in this process. It follows that
$$
b_n = \sum_{\ell \geq 1} \sum\left\{\prod_{i=1}^\ell a_i : \mbox{compositions $k_1+\ldots + k_\ell=n$}\right\},
$$  
and so, setting $B(x)=\sum_{k \geq 0} b_kx^k$ and $A(x)=\sum_{\ell \geq 1} a_\ell x^\ell$, we have 
\begin{equation} \label{AB}
B(x) = \frac{1}{1-A(x)}.
\end{equation}
Notice that $A(x)=1$ has a unique solution $r_{\rm pos}$ in the interval $(0,1)$. 
Let $r$ be any upper bound on $r_{\rm pos}$. From (\ref{AB}), standard facts about generating functions (see e.g. \cite[Section 2.4]{Wilf}) tells us that $r_{\rm pos}$ (and thus $r$) is an upper bound on the radius of convergence of $B(x)$; consequently
$$
\limsup_{n \rightarrow \infty} b_n^{1/n} \geq 1/r.  
$$
But we also know that $\c_n \ge b_n$, implying that 
$$
\mu_{\rm taxi} = \limsup_{n \rightarrow \infty} \c_n^{1/n} \ge \limsup_{n \rightarrow \infty} b_n^{1/n} \ge 1/r. 
$$
It follows that an upper bound on $r_{\rm pos}$, the unique positive solution to $A(x)=1$, yields a lower bound on $\mu_{\rm taxi}$.

Consider a sequence $(a_n')_{n=1}^\infty$ with $0 \leq a_n' \leq a_n$ for each $n$. Set $A'(x)=\sum_{n=1}^\infty a_n'x^n$. As is the case with $A(x)$, the equation $A'(x)=1$ has a unique solution $r'$ in the interval $(0,1)$, which moreover clearly satisfies $r' \geq r_{\rm pos}$. We record the conclusion to this discussion as a theorem.
\begin{theorem} \label{thm-using_irr_bridges}
With the notation as above, if $x>0$ satisfies $\sum_{n=1}^\infty a_n' x^n >1$ then $\mu_{\rm taxi} > 1/x$.
\end{theorem}

For any $N\geq 1$, the coefficients of the power series of $1/B(x)$ up to the coefficient of $x^N$ are determined by the coefficients of the power series of $B(x)$ up to $x^N$, and so using (\ref{AB}) the coefficients of $A(x)$ up to $x^N$ are determined by the coefficients of $B(x)$ up to $x^N$. 
We know the coefficients of $B(x)$ up to $x^{60}$ and from this we can easily calculate $a_n$ for $n \leq 60$ (see \cite{BGRT-data}). 
Taking $a_n'=a_n$ for $n \leq 60$ and $a_n'=0$ for $n > 60$ yields $\mu_{\rm taxi} > 1.55701$ and $\mu_{\rm taxi}^4-1 > 4.8771$, as stated in Theorem \ref{thm-best_mut}. 

\section{Concluding remarks} \label{sec-remarks}

\begin{itemize}

\item In an early version of this work \cite{BGRT} we employed the Peierls argument described in Section \ref{sec-proof_phase_coexistence} to separate $\mu_n^{\rm e}$ and $\mu_n^{\rm o}$ by establishing  
$$
\mu_n^{\rm e}(E) < 1/3
$$
in place of (\ref{distinguishing-statistic-0}), where the event $E$ was defined in terms of fault lines and crosses, as defined in \cite{randall}. The idea of replacing $E$ with ${\mathcal O}_m$, and conditioning on ${\mathcal H}_m$, which significantly streamlines the analysis, was suggested to us by Koteck\'y \cite{kot-pc} after the third author spoke on this work during the 2013-14 Warwick EPSRC Symposium on Statistical Mechanics. We are very grateful to him for this suggestion.

\item The standard Peierls argument for establishing phase coexistence tries to separate $\mu_n^{\rm e}$ and $\mu_n^{\rm o}$ using the event $\{v \in I\}$ where $v \in {\mathcal O}$ is some fixed vertex. Indeed, in \cite{bs} it is shown using the FKG inequality that there is phase coexistence for the hard-core model on ${\mathbb Z}^2$ if {\em and only if} $\limsup_{n \rightarrow \infty} \mu_n^{\rm e}(\{v \in I\}) < \limsup_{n \rightarrow \infty} \mu_n^{\rm o}(\{v \in I\})$. However, analyzing $\mu_n^{\rm e}(\{v \in I\})$ using the approach described in the present paper requires considering contours of length $4\ell$ for all $\ell \geq 3$. This in turn necessitates controlling the sub-exponential term in the growth rate of taxi walks, which in turn leads to poorer bounds on $\lambda$. Using ${\mathcal O}_m$ as the distinguishing event, which ensures that all contours are long, obviates this necessity.
This specifically comes in to play with the lower bound on $\ell$ in (\ref{long_contours_X}), coming from Lemma \ref{lem-contours_are_long4}.    

\item In \cite{BGRT} it was shown that if $\lambda$ satisfies $\lambda > \mu_{\rm taxi}^4-1$ then Glauber dynamics for independent sets
on the $n \times n$ torus takes time at least $e^{cn}$ to mix, for some constant $c = c(\lambda) > 0$, and that if also $2(1 + \lambda) > \mu_{\rm taxi}^2(1 + \sqrt{1 + 4\lambda})$ then the same is true for the $n \times n$ grid. Based on the computations in that reference, it was possible to conclude slow mixing for $\lambda>5.3646$ on the torus, and $\lambda >  7.1031$ on the grid. Using our improved bounds on $\mu_{\rm taxi}$ here, we can improve these bounds to $\lambda>\besttorus$ on the torus, and $\lambda >  \bestbox$ on the grid.  

\item One way to improve our lower bound on $\mu_{\rm taxi}$ would be to construct families of irreducible bridges of various lengths $n > 60$, and to use the sizes of these families as the $a_n'$s in Theorem \ref{thm-using_irr_bridges}. So far we have only had slight success with this approach, obtaining $\mu_{\rm taxi} > 1.55711$ and $\mu_{\rm taxi}^4-1> 4.8786$. The details are messy, and we choose not to include them here, but they can be found at \cite{BGRT-data}. 

\end{itemize}

\end{document}